\documentclass[11pt]{article}
\usepackage{amsmath,amsthm,amsfonts,amscd,bezier}
\usepackage{amssymb}
\usepackage{enumerate}
\usepackage[T1]{fontenc}
\usepackage[usenames]{color}

\newtheorem{lemma}{Lemma}[section]

\newtheorem{theorem}[lemma]{Theorem}
\newtheorem{remark}[lemma]{Remark}
\newtheorem{proposition}[lemma]{Proposition}

\newtheorem{definition}[lemma]{Definition}
\newtheorem{example}[lemma]{Example}
\newcommand{\Dem}{\noindent{\sc Proof:\ \ }}
\newcommand{\cqd}{{\hfill $\rule{2mm}{2mm}$}\vspace{1cm}}

\title{On the value set of $1$-forms for plane branches}

\author{{\sc Marcelo Osnar Rodrigues de Abreu and Marcelo Escudeiro Hernandes}\thanks{The first author was partially supported by CAPES and the second one by CNPq-Brazil.} \thanks{Corresponding author: Hernandes, M. E.; email: mehernandes@uem.br }}

\begin{document}
\maketitle \markboth{M. O. R. Abreu and M. E. Hernandes}{The value
set of $1$-forms determine the value semigroup for plane curves}

\begin{center} 2010 Mathematics Subject Classification: 14H20 (primary),
 32S10 (secondary)\end{center}

\begin{center} keywords: Plane curves, Semigroup, Value set of $1$-forms. \end{center}

\begin{abstract}
The value semigroup $\Gamma$ and the value set $\Lambda$ of
$1$-forms are, respectively, a topological and an analytical
invariant of a plane branch. Giving a plane branch $\mathcal{C}$ with
semigroup $\Gamma$ there are a finitely number of distinct possible
sets $\Lambda_i$ according to the analytic class of $\mathcal{C}$. In this
work we show that the value set of $1$-forms
$\Lambda$ determines the semigroup $\Gamma$ and we present an
effective method to recover $\Gamma$ by $\Lambda$. In particular,
this allows us to decide if a subset of $\mathbb{N}$ is a value set of
$1$-forms for a plane branch.
\end{abstract}

\section{Introduction}

For a plane branch $\mathcal{C}$ given by $f=0$ with
$f\in\mathbb{C}\{x,y\}$ irreducible, the topological class is totaly
characterized by numerical invariants, for instance, by the
characteristic exponents of a Puiseux parametrization
$\varphi(t)\in\mathbb{C}\{t\}\times\mathbb{C}\{t\}$
or, equivalently, by the value semigroup
$\Gamma_{\mathcal{C}}=\{\nu_f(h):=ord_t \varphi^*(h);\
h\in\mathbb{C}\{x,y\}\setminus\langle f\rangle\}\subseteq\mathbb{N}$, where $\varphi^*(h)=h(\varphi)$
(see \cite{equi}).

On the other hand, considering
$\Omega^1=\mathbb{C}\{x,y\}dx+\mathbb{C}\{x,y\}dy$ the
$\mathbb{C}\{x,y\}$-module of holomorphic $1$-forms on $\mathbb{C}^2$ we
have that
$\Lambda_{\mathcal{C}}=\{\nu_f(\omega):=ord_t(t\cdot \varphi^*(\omega));\
\varphi^*(\omega)\neq 0\}$ is an analytic invariant of $\mathcal{C}$
and
$\Gamma_{\mathcal{C}}\setminus\{0\}\subseteq\Lambda_{\mathcal{C}}$.
Some classical results
about holomorphic $1$-forms on a complete intersection curve
singularities and relations with other analytic invariants are summarized by
Dimca and Greuel in \cite{dimca}.

In \cite{algoritmo}, we find an algorithm to compute the set $\Lambda$ for any irreducible curve plane or not. Such algorithm allows us to obtain all the
possible set $\Lambda$ for plane branches in a given topological
class determined by a value semigroup $\Gamma$ obtaining precisely the algebraic relations among the coefficients of a Puiseux parametrization in order to achieve a set $\Lambda$ and it allows us describe parameterizations for such plane branches. In \cite{delorme}, Delorme consider such subject for plane branches with semigroup $\Gamma_{\mathcal{C}}=\langle v_0,v_1\rangle$. In a more general situation,
Alberich, Almir\'on and Moyano-Fern\'andez in \cite{patricio} present a characterization for the value set of $\mathcal{O}$-modules $\mathcal{O}+z\mathcal{O}$ where $\mathcal{O}$ is the local ring of a plane branch $\mathcal{C}$ with semigroup $\Gamma_{\mathcal{C}}=\langle v_0,v_1\rangle$.

Many authors have
considered simultaneously the sets $\Gamma$ and $\Lambda$ to
establish results. For instance, Carbonne in \cite{carbonne}
introduce the equidifferentiability class of plane branches: two plane branches $\mathcal{C}$ and $\mathcal{D}$ are equidifferentiable if $\Gamma_{\mathcal{C}}=\Gamma_{\mathcal{D}}$ and $\Lambda_{\mathcal{C}}=\Lambda_{\mathcal{D}}$. In \cite{HefezHern} and \cite{handbook} is
presented normal forms with respect to the analytic equivalence for
plane branches with $\Gamma$ and $\Lambda$ fixed.

Plane branches with same value semigroup can admit distinct value
set of $1$-form. In this paper we show that the numerical set
$\Lambda_{\mathcal{C}}$ allows us to recover the value semigroup of
$\mathcal{C}$ (see Theorem \ref{temosgamma}) that is false for space
curves (see Example \ref{espacial})
and for plane curves with several branches (see Example
\ref{reduzida}).
As a consequence, we obtain a method to decide if a subset
$L\subset\mathbb{N}$ is a value set of $1$-forms of a plane branch,
that is, if $L=\Lambda_{\mathcal{C}}$ for some plane branch
$\mathcal{C}$ (see Algorithm 2).

Notice that, in particular, the notion of equidifferentiability of two plane branches $\mathcal{C}$ and $\mathcal{D}$, by Carbonne, can be achieved just by the equality $\Lambda_{\mathcal{C}}=\Lambda_{\mathcal{D}}$ and to obtain the analytic normal form in \cite{HefezHern} it is sufficient to fix the set $\Lambda_{\mathcal{C}}$.

\section{Semigroup and value set of $1$-forms}

In this paper, we consider $\mathbb{N}=\mathbb{Z}_{\geq 0}$ and
$\mathbb{C}\{x,y\}$ denotes the ring of analytic power series in the indeterminates $x$
and $y$ over the complex field $\mathbb{C}$.

Let $\mathcal{C}$ be (a germ of) a plane branch, i.e. an irreducible plane curve, defined by $f=0$ with $f\in\mathcal{M}:=\langle x,y\rangle\subset\mathbb{C}\{x,y\}$ irreducible. We refer to \cite{zariskibook} and \cite{abramo} for details to this section.

We say that two plane branches $\mathcal{C}$ and $\mathcal{D}$ are
topologically equivalent if there exist neighborhoods $V, U\subseteq
(\mathbb{C}^2,0)$ and a homeomorphism $\Phi
:(\mathbb{C}^2,0)\rightarrow (\mathbb{C}^2,0)$ such that
$\Phi(\mathcal{C}\cap V)=\mathcal{D}\cap U$. In addition, if $\Phi$
is an analytic isomorphism then we say that $\mathcal{C}$ and
$\mathcal{D}$ are analytically equivalent. We can rewrite the analytic equivalence by $\mathbb{C}\{x,y\}$-automorphism. Two curves given by $f=0$
and $g=0$ respectively, are analytically equivalent if and only if
there exist an automorphism $\Psi$ and a unit $u$ of
$\mathbb{C}\{x,y\}$ such that $\Psi(f)=u\cdot g$, or equivalent
$\mathcal{O}_f:=\frac{\mathbb{C}\{x,y\}}{\langle
f\rangle}\simeq\frac{\mathbb{C}\{x,y\}}{\langle
g\rangle}=:\mathcal{O}_g$ as $\mathbb{C}$-algebras.

Up to a linear automorphism $\Psi\in Aut(\mathbb{C}\{x,y\})$ and, by
the Weierstrass Theorem Preparation, any plane curve is analytically
equivalent to a curve defined by
$f=y^n+\sum_{i=1}^{n}c_i(x)y^{n-i}\in\mathbb{C}\{x\}[y]$ where
$n:=mult(f)$, that is,
$f\in\mathcal{M}^n\setminus\mathcal{M}^{n+1}$. The Newton-Puiseux
theorem allows us to determine the roots of $f$. More explicitly, if
$f\in\mathbb{C}\{x\}[y]$ is irreducible then we get $y=\sum_{i\geq
n}a_ix^{\frac{i}{n}}\in\mathbb{C}\left \{x^{\frac{1}{n}}\right \}$
such that $f=\prod_{j=1}^{n}\left (y-\sum_{i\geq
n}a_i(\zeta^jx^{\frac{1}{n}})^i\right )$ where $\zeta\in\mathbb{C}$ is a
primitive $n$-th root of unity. We call $\varphi(t):=\left
(t^n,\sum_{i\geq n}a_it^i\right )$ a Puiseux parametrization of the
plane branch defined by $f$.

The analytic equivalence of plane branches can also be translate in terms
of the $\mathcal{A}$-equivalence of their parametrization, that is,
$\mathcal{C}$ and $\mathcal{D}$ are analytical equivalent plane
branches with parametrizations $\varphi(t)$ and $\psi(t)$
respectively, if and only if there exist diffeomorphisms $\sigma\in
\mbox{Diff}(\mathbb{C}^2,0)$ and $\rho\in \mbox{Diff}(\mathbb{C},0)$ such that
$\sigma\circ\varphi\circ\rho^{-1}(t)=\psi(t)$.

Given a Puiseux parametrization $\varphi(t)=(t^n,\sum_{i\geq n}a_it^i)$ for a plane branch $\mathcal{C}$ we define
$$\beta_0:=n=:e_0,\ \ \ \beta_i:=\min\{j;\ a_j\neq 0\ \mbox{and}\ e_{i-1}\nmid j\}\ \ \ \mbox{and}\ \ \ e_i=gcd(e_{i-1},\beta_i)\ \ \ \mbox{for}\ \ \ i\geq 1.$$
There exists $g\geq 1$ such that $e_g=1$ then we get $\beta_0<\beta_1<\cdots <\beta_g$ and $1=e_g<\cdots <e_0$. The sequence $(\beta_i)_{i=0}^{g}$ is called
 the characteristic sequence and Zariski, in \cite{equi}, shows that it determines and it is determined by the topological class of $\mathcal{C}$.

The topological class of a plane branch $\mathcal{C}$ defined by $f\in\mathbb{C}\{x,y\}$ can also be characterized by others numerical data, for instance by the value semigroup. More explicitly, if $\varphi(t)=(t^n,y(t))$ is a Puiseux parametrization  of $\mathcal{C}$ then we have the exact sequence
$$\begin{array}{ccccccccc}
\{0\} & \rightarrow & \langle f\rangle & \rightarrow & \mathbb{C}\{x,y\} & \stackrel{\varphi^*}{\rightarrow} & \mathbb{C}\{t^n,y(t)\} & \rightarrow & \{0\} \\
& & & & h & \mapsto & \varphi^*(h):=h(t^n,y(t)) & & \end{array}$$
and
$\mathcal{O}_f\approx\mathbb{C}\{t^n,y(t)\}\subseteq\mathbb{C}\{t\}$.
The value semigroup of $\mathcal{C}$ is
$$\Gamma_{\mathcal{C}}=\left \{\nu_f(h):=ord_t(\varphi^*(h))=
\dim_{\mathbb{C}}\frac{\mathbb{C}\{x,y\} }{\langle f,h\rangle};\
h\in \mathbb{C}\{x,y\}\setminus\langle f\rangle\right
\}\subseteq\mathbb{N}.$$

In fact, $\Gamma_{\mathcal{C}}$ is a numerical semigroup with minimal set of generators $\{v_0,v_1,\ldots ,v_g\}$ satisfying
\begin{equation}\label{vi}v_0:=\beta_0,\ \ \ v_{i}:=n_{i-1}v_{i-1}+\beta_{i}-\beta_{i-1},\ \ \ \mbox{where}\ \ \ n_0:=1,\ n_i:=\frac{e_{i-1}}{e_i}\ \ \ \mbox{for}\ \ \ 1\leq i\leq g.\end{equation}

Remark that $e_i=gcd(v_0,\ldots ,v_i)$ for $0\leq i\leq g$ we have that
$\Gamma_{\mathcal{C}}=\langle v_0,\ldots ,v_g\rangle$ and the
characteristic sequence $(\beta_i)_{i=0}^{g}$ are mutually
determined (see \cite{zariskibook}).

A monic polynomial $f_k\in\mathbb{C}\{x\}[y]$ of degree
$n_0\cdot\ldots\cdot n_k=\frac{v_0}{e_k}$ with $0\leq k<g$ is called
a $k$-semiroot of $f$ if $\nu_f(f_k)=v_{k+1}$. It follows that $f_k$
is irreducible and the value semigroup associated to
$\mathcal{C}_{f_i}$ is $\Gamma_i=\langle \frac{v_0}{e_i},\ldots
,\frac{v_i}{e_i}\rangle$ (see \cite{popescu}).

Any set $\{\varphi^*(h_i);\ h_i\in\mathbb{C}\{x,y\},\ i=0,\ldots
,g\}\subseteq \mathbb{C}\{t\}$ such that $\nu_f(h_i)=v_i$ is called a
(minimal) Standard Basis for $\mathcal{O}_f$. We can obtain a
minimal Standard Basis for $\mathcal{O}_f$ considering $\{x,\ f_k,\
0\leq k<g;\ f_k\ \mbox{is a semiroot of}\ f\}$, by minimal
polynomial of particular truncate of $\varphi(t)$ (see
\cite{zariskibook}) or by an algorithm presented in
\cite{algoritmo}, for instance.

Given a value semigroup $\Gamma_{\mathcal{C}}=\langle v_0,\ldots
,v_g\rangle$ any integer $z\in\mathbb{Z}$ can be uniquely expressed
as $z=\sum_{i=0}^{g}s_iv_i$ with $s_0\in\mathbb{Z}$ and $0\leq
s_i<n_i$ for $1\leq i\leq g$. Moreover, $z\in\Gamma_{\mathcal{C}}$
if and only if $s_0\in\mathbb{N}$ (see Lemma 7.1, \cite{abramo}). In
particular,
$\mu-1:=\sum_{i=1}^{g}(n_i-1)v_i-v_0\not\in\Gamma_{\mathcal{C}}$ and
$\mu+\mathbb{N}\subseteq\Gamma_{\mathcal{C}}$. We call $\mu$ the
conductor of $\Gamma_{\mathcal{C}}$.

Notice that, by the above properties, we have that
$z\in\Gamma_{\mathcal{C}}$ if and only if
$\mu-1-z\not\in\Gamma_{\mathcal{C}}$, that is,
$\Gamma_{\mathcal{C}}$ is a symmetric semigroup.

The following theorem characterizes the semigroup of a plane branch.

\begin{theorem}[Bresinsky]\label{bresinsky} Given $0<v_0<\ldots < v_g$ with $gcd(v_0,\ldots ,v_g)=1$
we put $n_0=1$ and $n_i=\frac{gcd(v_0,\ldots
,v_{i-1})}{gcd(v_0,\ldots ,v_i)}$ for $1\leq i\leq g$. The semigroup
$\langle v_0,\ldots ,v_g\rangle$ is a value semigroup of a plane
branch if and only if $n_i\geq 2$ and $n_{i-1}v_{i-1}<v_{i}$ for all
$1\leq i\leq g$.
\end{theorem}
\Dem See Theorem 2 in \cite{bresinsky}. \cqd

On the other hand analytic invariants for plane branches are not
easy to compute. Some of them are related with the
$\mathbb{C}\{x,y\}$-module
$\Omega^1=\mathbb{C}\{x,y\}dx+\mathbb{C}\{x,y\}dy$.

Given a parametrization $\varphi(t)=(x(t),y(t))$ of a plane branch
$\mathcal{C}$ and $\omega=A(x,y)dx+B(x,y)dy\in\Omega^1$ we define
$\varphi^*(\omega):=\varphi^*(A)x'(t)+\varphi^*(B)y'(t)$ where $x'(t)$ and $y'(t)$ denote the derivative of $x(t)$ and $y(t)$ respectively. In this
way, we get the value set of $1$-forms
$$\Lambda_{\mathcal{C}}=\{\nu_f(\omega):=ord_t(t\cdot\varphi^*(\omega));\ \omega\in\Omega^1\ \ \mbox{and}\ \
\varphi^*(\omega)\neq 0\}\subset\mathbb{N}$$ that is an analytic
invariant.

As $\nu_f(dh)=\nu_f(h)$ for any $h\in\mathcal{M}$ it follows that
$\Gamma^{*}_{\mathcal{C}}:=\Gamma_{\mathcal{C}}\setminus\{0\}\subseteq\Lambda_{\mathcal{C}}$
and it is immediate that
$\Gamma_{\mathcal{C}}+\Lambda_{\mathcal{C}}\subseteq\Lambda_{\mathcal{C}}$,
that is, $\Lambda_{\mathcal{C}}$ is a
$\Gamma_{\mathcal{C}}$-monomodule.

Carbonne, in \cite{carbonne}, say that two plane branches
$\mathcal{C}$ and $\mathcal{D}$ are equidifferentiable if
$\Gamma_{C}=\Gamma_{D}$ and $\Lambda_{C}=\Lambda_{D}$. Two
analytically equivalent plane curves are equidifferentiability and
plane branches equidifferentiable are obviously topologically
equivalent, but the converse are clearly false. In fact, we have the
following result:

\begin{theorem}[Hefez, Hernandes]\label{normal} Let $\mathcal{C}$ be a plane branch with
value semigroup $\Gamma=\langle v_0,\ldots ,v_g\rangle$ and set of
values of $1$-forms $\Lambda$. If $\Lambda\setminus\Gamma\neq
\emptyset$, then $\mathcal{C}$ is analytically equivalent to a plane
branch with parametrization
$$
\bigg ( t^{v_0}, t^{v_1}+a_\lambda
t^{\lambda}+\sum_{{i>\lambda}\atop{i\not\in\Lambda-v_0}}a_it^i\bigg
),
$$
where $\lambda=\min(\Lambda\setminus \Gamma)-v_0$. Otherwise
$\mathcal{C}$ is analytically equivalent to a branch with
parametrization $(t^{v_0},t^{v_1})$. Moreover, a branch
parameterized by $( t^{v_0}, t^{v_1}+b_\lambda
t^{\lambda}+\sum_{\lambda< i\not\in\Lambda-v_0}b_it^i)$ is
analytically equivalent to $\mathcal{C}$ if, and only if, there
exists $\alpha\in \mathbb{C}^*$ such that $b_i=\alpha^{i-v_1}a_i$,
for all $i$.
\end{theorem}
\Dem See Theorem 2.1 in \cite{HefezHern} or Theorem 1.3.17 in
\cite{handbook}. \cqd

\begin{remark}\label{lambda} Notice that if
$\Lambda_{\mathcal{C}}\setminus\Gamma_{\mathcal{C}}\neq\emptyset$
then
$\min(\Lambda_{\mathcal{C}}\setminus\Gamma_{\mathcal{C}})>v_1+v_0$.
\end{remark}

We can compute the set $\Lambda_{\mathcal{C}}$ associated to a plane
branch $\mathcal{C}$ given by a parametrization $\varphi(t)$ using
the Algorithm $4.10$ presented in \cite{algoritmo}. For commodity to
the reader we summarize the main ingredients of the mentioned
algorithm.

A minimal Standard Basis for $\varphi^*(\Omega^1)$ is a set
$G=\{\varphi^*(\omega_i);\ \omega_i\in\Omega^1,\ 1\leq i\leq
s\}\subset\mathbb{C}\{t\}$ such that
$\Lambda_{\mathcal{C}}=\bigcup_{i=1}^{s}\left
(\varphi^*(\omega_i)+\Gamma_{\mathcal{C}}\right )$ and
$\varphi^*(\omega_i)\not\in
\varphi^*(\omega_j)+\Gamma_{\mathcal{C}}$ for $i\neq j$.

Let $\{\varphi^*(h_0),\ldots ,\varphi^*(h_g)\}$ be a minimal
Standard Basis for $\mathcal{O}_f$ and $\emptyset \neq K\subset
\varphi^*(\Omega^1)$, we say that $\varphi^*(\omega_r)$ is a
reduction of $\varphi^*(\omega)\in \varphi^*(\Omega^1)$ modulo $K$
if there exist $a\in \mathbb{C}$, $\alpha_0,\ldots
,\alpha_g\in\mathbb{N}$ and $\varphi^*(\omega_k)\in K$ such that
$$\omega_r=\omega-a\prod_{i=0}^{g}h_i^{\alpha_i}\omega_k,\ \
\mbox{with}\ \  \nu_f(\omega_r)>\nu_f(\omega)\ \ \mbox{or}\ \
\omega_r=0.$$ We write $\varphi^*(\omega) \;
\stackrel{K}{\longrightarrow} \; \varphi^*(\omega_r)$ if
$\varphi^*(\omega_r)$ is a final reduction of $\varphi^*(\omega)$
modulo $K$, that is, $\varphi^*(\omega_r)$ is obtained from
$\varphi^*(\omega)$ through a chain (possibly infinite) of
reductions, modulo $K$, and cannot be reduced further.

A minimal $S$-process of a pair of elements
$\varphi^*(\omega_p),\varphi^*(\omega_q)\in \varphi^*(\Omega^1)$ is
$$\varphi^*\left ( S(\omega_p,\omega_q)\right )=\varphi^*\left (
a\prod_{i=0}^{g}h_i^{\alpha_i}\omega_p+b\prod_{i=0}^{g}h_i^{\gamma_i}\omega_q\right
),$$ where $a,b\in \mathbb{C},$
$\nu_f(S(\omega_p,\omega_q))>\nu_f(\prod_{i=0}^{g}h_i^{\alpha_i}\omega_p
)=\nu_f(\prod_{i=0}^{g}h_i^{\gamma_i}\omega_q)$ and $(\alpha_0,\ldots
,\alpha_{g},\gamma_0,\ldots ,\gamma_{g})$ is a minimal solution of
the linear Diophantine equation
$\sum_{i=0}^{g}v_i\alpha_i+\nu_f(\omega_q)=\sum_{i=0}^{g}v_i\gamma_i+\nu_f(\omega_p)$.

The following algorithm gives us Standard Bases for
$\varphi^*(\Omega^1)$ (see Algorithm $4.10$ in \cite{algoritmo}):

\begin{center}
\noindent {\bf Algorithm 1:} Standard basis for
$\varphi^*(\Omega^1)$:

\begin{tabular}{|l|}
\hline \texttt{Input}: $H=\{h_j,\ j=0,\dots ,g\}$ a minimal Standard Basis for $\mathcal{O}_f$; \\
\texttt{Define:} $G_0=\emptyset$; $G_1:=\{\varphi^*(dh_j),\ j=0,\dots ,g\}$
and $i:=1;$ \\
\texttt{While} $G_i\neq G_{i-1}$ \texttt{do} \\
\ \ \ \ \texttt{Compute} $S:= \{ s;\ s$ is a minimal $S$-process of $G_i \};$\\
\ \ \ \ \ \ \ \ \ \ \ \ \ \ \ $R:= \{ r;\ s \; \stackrel{{G_i}}{\longrightarrow}\; r$ and
$r\neq 0,\; \forall s\in S \};$ \\
\ \ \ \ \texttt{Define:} $G_{i+1}:=G_i\cup R;$ \\
\ \ \ \ \ \ \ \ \ \ \ \ \ \ \ \ $i:=i+1;$ \\
\texttt{Output}: $G=\cup_{i\geq 0}G_i.$ \\
\hline
\end{tabular}
\end{center}

From any Standard Basis of $\varphi^*(\Omega^1)$ we can extract a
minimal Standard Basis and consequently, we obtain the set
$\Lambda_{\mathcal{C}}$. In addition, we can apply the above
algorithm to the family of all plane branches with a fixed semigroup
and consequently we obtain all possible value set of $1$-forms for
plane branches in a fixed topological class.

\begin{example}[\cite{6919}]\label{ex6919}
Let $\Gamma=\langle 6,9,19\rangle$ be a numerical semigroup, by
Theorem \ref{bresinsky} it is a value semigroup of plane branches.
By (\ref{vi}) we get the characteristic sequence $\beta_0=6,
\beta_1=9$ and $\beta_2=10$ and consequently any plane branch with
semigroup $\Gamma$ can be done by a parametrization
$(t^6,t^9+\sum_{i\geq 10}a_it^i)$ with $a_{10}\neq 0$. In
particular, $ord_t(t\cdot\varphi^*(2xdy-3ydx))=10=\min(\Lambda\setminus\Gamma)-6$.

As $\Gamma\setminus\{0\}\subset\Lambda$ for any possible set
$\Lambda$, by Theorem \ref{normal}, any plane branch in this
topological class is analytically equivalent to a branch with
parametrization
$$(t^6,t^9+t^{10}+a_{11}t^{11}+a_{14}t^{14}+a_{17}t^{17}+a_{20}t^{20}+a_{23}t^{23}+a_{26}t^{26}).$$
Applying the Algorithm 1 we obtain:
\begin{center}
    \begin{tabular}{|c|c|}
    \hline
    restrictions & $\Lambda\setminus\Gamma$ \\
    \hline
$a_{11}\not\in\{-\frac{1}{2},\frac{29}{18}\}$ & $\{16, 22, 26, 29, 32, 35, 41\}$ \\
\hline
$a_{11}=\frac{29}{18}$ & $\{16, 22, 26, 32, 35, 41\}$ \\
\hline
$a_{11}=-\frac{1}{2},\ \ 1152a^2_{14}-769a_{14}+1064a_{17}\neq 28$ & $\{16, 22, 29, 32, 35, 41\}$ \\
\hline
$a_{11}=-\frac{1}{2},\ \ 1152a^2_{14}-769a_{14}+1064a_{17}= 28$ & $\{16, 22, 29, 35, 41\}$  \\
\hline
    \end{tabular}
\end{center}
\end{example}

\section{The characterization of the set $\Lambda$ for plane branches}

The aim of this section is to show that the value set
$\Lambda_{\mathcal{C}}$ of $1$-forms on a plane branch $\mathcal{C}$
determines the value semigroup $\Gamma_{\mathcal{C}}$ and to present
a method to decide is a subset $L\subset\mathbb{N}$ is a value set
of $1$-forms for some plane branch. As
$\min(\Lambda_{\mathcal{C}})=\min(\Gamma^{*}_{\mathcal{C}})=v_0$ we have that $v_0=1$ if and only if
$\Lambda_{\mathcal{C}}=\Gamma^{*}_{\mathcal{C}}=\mathbb{N}\setminus\{0\}$.
In what follows we assume that $v_0>1$.

First of all we present some arithmetic notions for a nonempty set
$S\subseteq\mathbb{N}\setminus\{0\}$.

Let $S$ be a nonempty subset of $\mathbb{N}\setminus\{0\}$. We say
that $S$ {\em admits a conductor} $s$ if $s-1\not\in S$ and
$s'\in S$ for every $s'\geq s$, or equivalently if
$\mathbb{N}\setminus S$ is a finite set.

\begin{definition}
Given $\emptyset\neq S\subseteq\mathbb{N}\setminus\{0\}$ we define
$$\mathbf{a}_0:=\min (S)\ \ \ \mbox{and}\ \ \ \mathbf{a}_i:=\min\left (S\setminus\bigcup_{j=0}^{i-1}(\mathbf{a}_{j}+\mathbb{N}\cdot \mathbf{a}_0)\right )\ \ \mbox{for}\ \ i>0.$$
As $\mathbf{a}_j\not\equiv \mathbf{a}_i\mod \mathbf{a}_0$ for $i\neq j$ there exists $m<\mathbf{a}_0$
such $S\setminus\bigcup_{i=0}^{m}(\mathbf{a}_i+\mathbb{N}\cdot
\mathbf{a}_0)=\emptyset$. The {\bf Ap\'{e}ry set} of $S$ is $Ap(S)=\{\mathbf{a}_i;\
0\leq i\leq m\}$.

If $Ap(S)=\{\mathbf{a}_0,\ldots ,\mathbf{a}_m\}$ then $m\leq \mathbf{a}_0-1$ and $S\subseteq
\bigcup_{j=0}^{m}(\mathbf{a}_j+\mathbb{N}\cdot \mathbf{a}_0)$. We say that $S$ is {\bf
covered by $Ap(S)$} if $m=\mathbf{a}_0-1$ and
$S=\bigcup_{j=0}^{m}(\mathbf{a}_j+\mathbb{N}\cdot \mathbf{a}_0)$.
\end{definition}

\begin{remark} If $S$ is covered by $Ap(S)$ and $\mathbf{a}_*=\max(Ap(S))$ then
$S$ has conductor $\mathbf{a}_*-\mathbf{a}_0+1$. In fact, by definition
$\mathbf{a}_*-\mathbf{a}_0\not\in S$ and for any $s>\mathbf{a}_*-\mathbf{a}_0$, as $\sharp Ap(S)=\mathbf{a}_0$
there exists $\mathbf{a}_j\in Ap(S)$ such that $\mathbf{a}_*-\mathbf{a}_0<s=\mathbf{a}_j+k \mathbf{a}_0\leq
\mathbf{a}_*+k \mathbf{a}_0$, consequently $k\in\mathbb{N}$ and $s\in
\mathbf{a}_j+\mathbb{N}\cdot \mathbf{a}_0\subseteq S$. The converse is false, for
instance $S=\mathbb{N}\setminus\{0,2\}$ admits conductor but it is
not covered by $Ap(S)$.
\end{remark}

\begin{remark}\label{classico} For any plane branch $\mathcal{C}$ the sets
$\Gamma^{*}_{\mathcal{C}}=\Gamma_{\mathcal{C}}\setminus\{0\}$ and
$\Lambda_{\mathcal{C}}$ are covered by their Ap\'{e}ry set. In fact, given $\Gamma_{\mathcal{C}}=\langle v_0,\ldots ,v_g\rangle$, by Proposition 7.11 in \cite{abramo}, we have that
$Ap(\Gamma^{*}_{\mathcal{C}})=\{v_0\}\ \dot{\cup}\
\{\sum_{i=1}^{g}s_iv_i\neq 0;\ 0\leq s_i<n_i\}$. Consequently,
$\sharp Ap(\Gamma^{*}_{\mathcal{C}})=\prod_{i=1}^{g}n_i=v_0$ and
$\Gamma_{\mathcal{C}}=\{\sum_{i=1}^{g}s_iv_i+\mathbb{N}\cdot
v_0;\ 0\leq s_i<n_i\}$.

As $\Gamma^{*}_{\mathcal{C}}\subseteq\Lambda_{\mathcal{C}}$ we have
that $Ap(\Lambda_{\mathcal{C}})=\{\mathbf{b}_i:=\mathbf{a}_i-k_iv_0;\ \mathbf{a}_i\in
Ap(\Gamma^{*}_{\mathcal{C}})\ \mbox{for some}\ k_i\in\mathbb{N}\}$, so
$\sharp Ap(\Lambda_{\mathcal{C}})=\sharp
Ap(\Gamma^{*}_{\mathcal{C}})=v_0$. Notice that
$\mathbf{b}_0:=\min(\Lambda_{\mathcal{C}})=\mathbf{a}_0=v_0$ and, by Remark
\ref{lambda}, we must have $\mathbf{b}_1:=\mathbf{a}_1=v_1$. As $\mathbb{N}\cdot
v_0+\Lambda_{\mathcal{C}}\subset\Gamma_{\mathcal{C}}+\Lambda_{\mathcal{C}}\subset\Lambda_{\mathcal{C}}$
it follows that
$\Lambda_{\mathcal{C}}=\bigcup_{j=0}^{v_0-1}(\mathbf{b}_j+\mathbb{N}\cdot
\mathbf{b}_0)$.
\end{remark}

In what follows, for any subset $S\subseteq\mathbb{N}\setminus\{0\}$
covered by its Ap\'{e}ry set $Ap(S)=\{\mathbf{a}_i;\ 0\leq i<\mathbf{a}_0\}$ we put
$\varepsilon_0:=\mathbf{a}_0, \eta_0:=1,$
\begin{equation}\label{para}\varepsilon_i:=\max\{gcd(\varepsilon_{i-1},\mathbf{a}_j); \ a_j\in
Ap(S)\ \mbox{and}\ \varepsilon_{i-1}\nmid \mathbf{a}_j\}\ \ \ \mbox{and}\ \ \
\eta_i:=\frac{\varepsilon_{i-1}}{\varepsilon_i}\ \ \ \mbox{for}\ \ \
i\geq 1.\end{equation}

Notice that there exists $\varrho\leq \mathbf{a}_0$
such that $\varepsilon_{\varrho}=1$, $\eta_i>1$ for any $1\leq
i\leq\varrho$ and $\prod_{j=0}^{i}\eta_j=\frac{\varepsilon_0}{\varepsilon_i}=\sharp\{\mathbf{a}\in Ap(S);\ \varepsilon_i\mid
\mathbf{a}\}$.

For $0<i\leq\varrho$ we denote $Q_1<\ldots <Q_{\frac{\varepsilon_0}{\varepsilon_{i-1}}}$ the $\frac{\varepsilon_0}{\varepsilon_{i-1}}$ smallest elements in
the set $\Delta_i=\{\mathbf{a}\in Ap(S);\ \varepsilon_i\mid \mathbf{a}\ \mbox{and}\
\varepsilon_{i-1}\nmid \mathbf{a}\}$ and we put
\begin{equation}\label{Bi} B_0(S):=\{\mathbf{a}_0\}\ \ \mbox{and}\ \  B_i(S):=\left \{Q_1,\ldots ,Q_{\frac{\varepsilon_0}{\varepsilon_{i-1}}}\right \}\ \ \mbox{for}\ \ 1\leq i\leq\varrho.\end{equation}

\begin{remark}\label{igual} For $S=\Gamma^{*}_{\mathcal{C}}$ or $S=\Lambda_{\mathcal{C}}$ we have
$\varepsilon_i=e_i$, $\eta_i=n_i$, $\varrho=g$, $B_0(S)=\{v_0\}$ and
$B_{1}(S)=\{v_1\}$. In particular, for $\varepsilon_1=e_1=1$ the
sets $B_0(\Lambda_{\mathcal{C}})$ and $B_1(\Lambda_{\mathcal{C}})$
determine $\Gamma_{\mathcal{C}}$.

In addition, by the description of $Ap(\Gamma^*_{\mathcal{C}})$
presented in Proposition 7.11 in \cite{abramo} (see Remark
\ref{classico} ) we have
$B_i(\Gamma^*_{\mathcal{C}})=\{v_i+\sum_{j=1}^{i-1}s_jv_j;\ 0\leq
s_j< n_j\}$ for $1\leq i\leq g$.
\end{remark}

By the above remark we get
$v_i=\min(B_i(\Gamma^{*}_{\mathcal{C}}))$ for $0\leq i\leq g$. In
that follows we will show that
$v_i=\max(B_i(\Lambda_{\mathcal{C}}))$ where $0\leq i\leq g$ for any set
of value of $1$-forms $\Lambda_{\mathcal{C}}$ associated to a plane
branch $\mathcal{C}$, consequently the set $\Lambda_{\mathcal{C}}$
determines the value semigroup $\Gamma_{\mathcal{C}}$.

From now on we will assume that $\Lambda_{\mathcal{C}}$ is a value
set of $1$-form for some plane branch $\mathcal{C}$ with value
semigroup $\Gamma_{\mathcal{C}}=\langle v_0,\ldots ,v_g\rangle$ with
$g\geq 2$.

\begin{lemma}\label{lema} There is no elements $v_{i}-kv_0$ in $\Lambda_{\mathcal{C}}$ with
$k>0$ and $0\leq i\leq g$.
\end{lemma}
\Dem First of all, notice that
$v_{i}-kv_0\not\in\Gamma_{\mathcal{C}}$ for any $k>0$. By Remark
\ref{lambda},
$\min(\Lambda_{\mathcal{C}}\setminus\Gamma_{\mathcal{C}})>v_1+v_0$
then we get the result for $i=0$ and $i=1$.

Consider $i\geq 2$, $f_{i-1}$ a $(i-1)$-semiroot with value
semigroup $\Gamma_{i-1}$ and value of $1$-forms $\Lambda_{i-1}$
associated to $\mathcal{C}_{f_{i-1}}$. Given
$\delta\in\Lambda_{i-1}\setminus\Gamma_{i-1}$ we denote
$\Theta_{i-1}(\delta)=\max\{\nu_{f_{i-1}}(B);\
\delta=\nu_{f_{i-1}}(Adx-Bdy)\}$,
$L^1_{i-1}=\{\delta\in\Lambda_{i-1}\setminus\Gamma_{i-1};\
e_{i-1}(\delta-\Theta_{i-1}(\delta))<\beta_{i}\}$ and
$L^2_{i-1}=\{\delta\in\Lambda_{i-1}\setminus\Gamma_{i-1};\
e_{i-1}(\delta-\Theta_{i-1}(\delta))>\beta_{i}\}$.

By Proposition 5.1 in \cite{osnar}, we have that
\begin{equation}\label{conjunto}
\{\lambda\in\Lambda_{\mathcal{C}}\setminus\Gamma_{\mathcal{C}};\
\lambda<v_i\}= \rho_{i-1}(\Lambda_{i-1}\setminus\Gamma_{i-1})\
\dot{\cup}\ \{v_i-e_{i-1}\delta;\
0\neq\delta\in\mathbb{N}\setminus\Lambda_{i-1}\},
\end{equation}
where (see Proposition 4.7 in \cite{osnar})
$$\rho_{i-1}(\Lambda_{i-1}\setminus\Gamma_{i-1})=e_{i-1}L^1_{i-1}\
\dot{\cup}\ \left \{v_i-e_{i-1}\left (
\mu_{i-1}-1+\frac{v_0}{e_{i-1}}-\Theta_{i-1}(\delta)\right );\ \delta\in
L^2_{i-1}\right \}$$ and $\mu_{i-1}$ is the conductor of
$\Gamma_{i-1}$.

Suppose that
$v_i-kv_0\in\rho_{i-1}(\Lambda_{i-1}\setminus\Gamma_{i-1})$. Notice
that $e_{i-1}\nmid v_i-kv_0$, consequently $v_i-kv_0\not\in
e_{i-1}L^1_{i-1}$. In addition, if $v_i-kv_0=v_i-e_{i-1}\left (
\mu_{i-1}-1+\frac{v_0}{e_{i-1}}-\Theta_{i-1}(\delta)\right )$ as
$\Theta_{i-1}(\delta)\in\Gamma_{i-1}$ and $\Gamma_{i-1}$ is a
symmetric semigroup we have
$(k-1)\frac{v_0}{e_{i-1}}=\mu_{i-1}-1-\Theta_{i-1}(\delta)\not\in\Gamma_{i-1}$
that is an absurd.

If $v_i-kv_0=v_i-e_{i-1}\delta$ then
$\delta=k\frac{v_0}{e_{i-1}}\in\Gamma_{i-1}\setminus\{0\}\subseteq\Lambda_{i-1}$,
so $v_i-kv_0\not\in\{v_i-e_{i-1}\delta;\
0\neq\delta\in\mathbb{N}\setminus\Lambda_{i-1}\}$.

Hence, $v_i-kv_0\not\in \Lambda_{\mathcal{C}}$ for $k>0$ and $0\leq
i\leq g$. \cqd

Now we can describe some properties of
$B_{i}(\Lambda_{\mathcal{C}})$ for $2\leq i\leq g$.

\begin{proposition} For $2\leq i\leq g$ we have that $v_i=\max (B_i(\Lambda_{\mathcal{C}}))$.
\end{proposition}
\Dem By Remark \ref{classico}, we get
$Ap(\Gamma^{*}_{\mathcal{C}})=\{v_0\}\ \dot{\cup}\
\{\sum_{i=1}^{g}s_iv_i\neq 0;\ 0\leq s_i<n_i\}$ and
$Ap(\Lambda_{\mathcal{C}})=\{\mathbf{b}_i:=\mathbf{a}_i-k_iv_0;\ \mathbf{a}_i\in
Ap(\Gamma^{*}_{\mathcal{C}})\ \mbox{for some}\ k\in\mathbb{N}\}$.

Denoting
$\Delta_i:=\{\alpha\in Ap(\Lambda_{\mathcal{C}});\ e_i\mid \alpha\ \mbox{and}\ e_{i-1}\nmid \alpha\}$, by the previous lemma we get $v_i\in \Delta_i$.

As $\min(\Lambda_{i-1})=\frac{v_0}{e_{i-1}}=n_0\cdots n_{i-1}$ it follows that $\left \{1,\ldots ,\frac{v_0}{e_{i-1}}-1\right \}\subset\mathbb{N}\setminus\Lambda_{i-1}$
and by (\ref{conjunto}) we get $v_i-e_{i-1}\delta\in\Lambda_{\mathcal{C}}\setminus\Gamma_{\mathcal{C}}$ for $1\leq \delta < \frac{v_0}{e_{i-1}}$. In addition, $v_i-e_{i-1}\delta_j\not\equiv v_i-e_{i-1}\delta_k\mod v_0$ for $1\leq \delta_k,\delta_j< n_0\cdots n_{i-1}$ and $\delta_j\neq\delta_k$. Consequently, $$\{v_{i}\}\ \cup\ \left \{v_i-e_{i-1}\delta-k_{\delta}v_0;\ \ 1\leq \delta<\frac{v_0}{e_{i-1}}\ \mbox{and some}\ k_{\delta}\in\mathbb{N}\right \}\subset \{\mathbf{a}\in\Delta_i;\ \mathbf{a}\leq v_i\}.$$

On the other hand, if $\alpha\in \Delta_i$ with $\alpha <v_i$ as $v_i=\min\{\gamma\in\Gamma_{\mathcal{C}};\ e_i\mid \gamma\ \mbox{and}\ e_{i-1}\nmid \gamma\}$ it follows that $\alpha\in\Lambda_{\mathcal{C}}\setminus\Gamma_{\mathcal{C}}$ and by (\ref{conjunto}) we get $\alpha=v_i-e_{i-1}\gamma$ with $\gamma>0$. As
$\alpha=v_i-e_{i-1}\gamma\equiv v_i-e_{i-1}\delta\mod v_0$ for some $1\leq \delta<\frac{v_0}{e_{i-1}}$ we get
$$\{\mathbf{a}\in\Delta_i;\ \mathbf{a}\leq v_i\}=\{v_{i}\}\ \cup\ \left \{v_i-e_{i-1}\delta-k_{\delta}v_0;\ \ 1\leq \delta<\frac{v_0}{e_{i-1}}\ \mbox{and some}\ k_{\delta}\in\mathbb{N}\right \}=B_i(\Lambda_{\mathcal{C}})$$
and hence $\max(B_i(\Lambda_{\mathcal{C}}))=v_i$ for all $2\leq i\leq g$.
\cqd

We can recover the semigroup $\Gamma_{\mathcal{C}}$ from the set $\Lambda_{\mathcal{C}}$ as we highlighted in the following result:

\begin{theorem}\label{temosgamma}
    The value set of $1$-forms $\Lambda_{\mathcal{C}}$ determines the value semigroup $\Gamma_{\mathcal{C}}$.
\end{theorem}
\Dem Given $\Lambda_{\mathcal{C}}$ we obtain $Ap(\Lambda_{\mathcal{C}})$ and we determine the sequences $(\varepsilon_i)_{i=0}^{g}$ and $(n_i)_{i=0}^{g}$ (see (\ref{para}) and Remark \ref{igual}). This allows us to determine $B_i(\Lambda_{\mathcal{C}})$ for $0\leq i\leq g$.

As $B_0(\Lambda_{\mathcal{C}})=\{v_0\}$, $B_1(\Lambda_{\mathcal{C}})=\{v_1\}$ and by the above proposition, we get  $v_i=\max(B_i(\Lambda_{\mathcal{C}}))$ for $0\leq i\leq g$, that is, the miminal generators of $\Gamma_{\mathcal{C}}$.
\cqd

As a consequence of the above arithmetic properties of the analytic invariant $\Lambda_{\mathcal{C}}$ we can conclude that if two plane branches have the same value set of $1$-form then the branches have the same value semigroup.
The following example illustrates that this property is false for space curves.

\begin{example}\label{espacial}
Consider the curves given by the parametrizations:
$$\mathcal{C}_1:\left \{
\begin{array}{l}
x=t^6\\
y=t^{14}+t^{17}\\
z=t^{39}
\end{array}
\right .\ \ \ \mbox{and}\ \ \ \mathcal{C}_2:\left \{
\begin{array}{l}
x=t^6\\
y=t^{14}+t^{33}\\
z=t^{23}.
\end{array}
\right .$$

The corresponding associated value semigroup, computed by Algorithm 3.2 in \cite{algoritmo} for instance, are
$\Gamma_1 =\langle 6,14,39\rangle$ and $\Gamma_2 =\langle 6,14,23\rangle$.

On the other hand, computing the value set of $1$-forms by Algorithm 4.10 in \cite{algoritmo}, we get $\nu(3xdy-7ydx)=
23$ for $\mathcal{C}_1$ and $\nu(3xdy-7ydx)=39$ for $\mathcal{C}_2$.

In this way, we obtain $\Lambda_1 =\Lambda_2 =\{0,6,12,14,18,20,23,24,26,28,29,30,32,34+\mathbb{N}\}$.
\end{example}

We can extend the value semigroup and the value set of $1$-forms for
plane curves with several branches. More explicitly, if
$\mathcal{D}=\cup_{i=1}^{r}\mathcal{C}_i$ and
$\varphi_i(t_i)=(x(t_i),y(t_i))$ denotes a parametrization of
$\mathcal{C}_i$ then
$$\Gamma_{\mathcal{D}}:=\{(ord_{t_1}(\varphi_1^*(h)),\ldots ,ord_{t_r}(\varphi_r^*(h)));\
 h\in\mathbb{C}\{x,y\}\setminus\cup_{i=1}^{r}\langle f_i\rangle \}\ \ \ \mbox{and}$$
$$\Lambda_{\mathcal{D}}:=\{(ord_{t_1}(t_1\cdot\varphi_1^*(\omega)),\ldots ,ord_{t_r}(t_r\cdot\varphi_r^*(\omega)));\
 \omega\in\Omega^1\ \mbox{and}\ \varphi_i^*(\omega)\neq 0,\ 1\leq i\leq r \}.$$
In this situation, the Theorem \ref{temosgamma} is not true for
$r>1$ as we can observe in the next example.

\begin{example}\label{reduzida}
Let $\mathcal{D}_{\alpha}=\mathcal{C}_1\cup\mathcal{C}_2$ be a plane
curve such that the branches are given by $f_1=y^2-2x^2y-x^3+x^4$
and $f_2=y^2-\alpha x^3$ with $\alpha\in\mathbb{C}\setminus\{0\}$,
or equivalently, given by the parametrizations
$\varphi_1(t)=(t_1^2,t_1^3+t_1^4)$ and
$\varphi_2(t)=(t_2^2,\sqrt{\alpha}\cdot t_2^3)$.

We can obtain $\Gamma_{\mathcal{D}_{\alpha}}$ and
$\Lambda_{\mathcal{D}_{\alpha}}$ by means algorithms in
\cite{carvalho}, but we choose to describe explicitly the computations.

Any element $h\in\mathbb{C}\{x,y\}$ can be expressed by
$h=q_1.f_1+q_2.f_2+A(x).y+B(x)$ and computing
$(ord_{t_1}(\varphi_1(h)),ord_{t_2}(\varphi_2(h))$ for any
$h\in\mathbb{C}\{x,y\}\setminus\langle f_1\rangle\cup\langle
f_2\rangle$ we get
$$\Gamma_{\mathcal{D}_{\alpha}}=\{(0,0), (a,a);\ 2\leq a\leq 7\}\cup \{(6,6+\mathbb{N}),
 (6+\mathbb{N},6), (8,8)+\mathbb{N}^2\} \ \ \mbox{for}\ \ 0\neq\alpha\neq 1\ \ \mbox{and}$$
$$\Gamma_{\mathcal{D}_{1}}=\{(0,0), (a,a);\ 2\leq a\leq 8\}\cup \{(7,7+\mathbb{N}),
 (7+\mathbb{N},7), (9,9)+\mathbb{N}^2\},$$
that is, $\Gamma_{\mathcal{D}_{\alpha}}\neq\Gamma_{\mathcal{D}_{1}}$
for $0\neq\alpha\neq 1$.

Now, analyzing the values of $1$-form we verify that the elements
$(a,b)\in\Lambda_{\mathcal{D}_{\alpha}}\setminus
\Gamma_{\mathcal{D}_{\alpha}}$ are such that $a,b\geq 6$.

Let us consider $\omega_1=3ydx-2xdy, \omega_2=\omega_1+x^2dx\in\Omega^1$.

For $0\neq\alpha\neq 1$ and for any $\gamma=2i+3j\in\langle
2,3\rangle$ the values of $df_2+3(1-\alpha)\omega_1+x^iy^jdf_1$ and
$df_1+3(1-\alpha)\omega_2+x^iy^jdf_2$ are $(7,6+\gamma)$ and
$(6+\gamma,7)$ respectively.

For $\alpha=1$ and for any $\gamma=2i+3j\in\langle 2,3\rangle$, the
$1$-forms $\omega_1+df_1, \omega_1+x^iy^j\omega_2,
\omega_2+df_2,\omega_2+x^iy^j\omega_1$, $x(\omega_1+df_1),
x(\omega_1+x^iy^j\omega_2), x(\omega_2+df_2)$ and
$x(\omega_2+x^iy^j\omega_1)$ give us the respective values $(6,7),
(6,6+\gamma), (7,6), (6+\gamma,6)$, $(8,9), (8,8+\gamma), (9,8)$ and
$(8+\gamma,8)$.

In this way, we have
$\Lambda_{\mathcal{D}_{\alpha}}=\{(2,2),(3,3),(4,4),(5,5),(6,6)+\mathbb{N}^2\}$
for any $\alpha\in\mathbb{C}\setminus\{0\}$.
\end{example}

Combining the previous result with Theorem \ref{bresinsky} and the
Algorithm 1 we have a method to verify if a subset
$L\subseteq\mathbb{N}\setminus\{0\}$ is a value set of $1$-forms for
some plane branch.

In fact, given $L\subseteq\mathbb{N}\setminus\{0\}$ we compute
$Ap(L)$, if $L$ is not covered by $Ap(L)$ then, by Remark
\ref{classico}, $L\neq\Lambda_{\mathcal{C}}$ for any plane branch
$\mathcal{C}$.

If $L$ is covered by its Ap\'{e}ry set, then we compute the
sequences $(\varepsilon_i)_{i=0}^{\varrho}$,
$(\eta_i)_{i=0}^{\varrho}$ (see (\ref{para})) and the sets $B_i(L)$ for $0\leq i\leq
\varrho$ (see (\ref{Bi})).

Taking $u_i=\max(B_i(L))$, if $(u_i)_{i=0}^{\varrho}$ do not
satisfy $\eta_i\geq 2$ and $\eta_{i-1}u_{i-1}<u_i$ for all $1\leq
i\leq \varrho$ then, by Theorem \ref{bresinsky}, the semigroup
$G:=\langle u_0,\ldots ,u_{\varrho}\rangle$ is not a value semigroup
for any plane branch $\mathcal{C}$ and, consequently by Theorem
\ref{temosgamma}, $L$ can not be a value set of $1$-forms for any
plane branch.

If the semigroup $G$ satisfies the conditions of Theorem
\ref{bresinsky}, that is, $G$ is a value semigroup for some
plane branch, then we compute all possible value sets $\Lambda_1,\ldots
,\Lambda_s$ of $1$-forms associated to it by the Algorithm $1$ and
the set $L$ will be a value set of $1$-forms for some plane branch
if and only if $L=\Lambda_j$ for some $1\leq j\leq s$.

We summarize the above discussion in the following algorithm.

\begin{center}
\noindent {\bf Algorithm 2:}

\begin{tabular}{|l|}
\hline \texttt{Input}: $L\subseteq\mathbb{N}\setminus\{0\}$; \\
\texttt{Compute:} $Ap(L);$ \\
\texttt{If} $L$ \texttt{is not covered by} $Ap(L)$ \\
\ \ \ \ \texttt{Then} $L$ is not a value set for any plane branch;\\
\ \ \ \ \texttt{Else} \\
\ \ \ \ \ \ \ \ \texttt{Compute:} $\eta_i,\ B_i(L)$
\texttt{and} $u_i=\max(B_i(L))$ \texttt{for} $0\leq i\leq \varrho;$ \\
\ \ \ \ \ \ \ \ \texttt{If} $\eta_i<2$ or $\eta_{i-1}u_{i-1}\geq u_i$ \texttt{for some} $1\leq i\leq\varrho$ \\
\ \ \ \ \ \ \ \ \ \ \ \ \texttt{Then}  $L$ is not a value set for any plane branch;\\
\ \ \ \ \ \ \ \ \ \ \ \ \texttt{Else} \\
\ \ \ \ \ \ \ \ \ \ \ \ \ \ \ \texttt{Use Algorithm 1 to compute all possible} $\Lambda_j$ \\
\ \ \ \ \ \ \ \ \ \ \ \ \ \ \ \texttt{for branches with semigroup} $\langle u_0,\ldots ,u_{\varrho}\rangle$; \\
\ \ \ \ \ \ \ \ \ \ \ \ \ \ \ \texttt{If} $L=\Lambda_j$ \texttt{for
some}
$j$\\
\ \ \ \ \ \ \ \ \ \ \ \ \ \ \ \ \ \ \ \texttt{Then} $L$ is a value set for some plane branch;\\
\ \ \ \ \ \ \ \ \ \ \ \ \ \ \ \ \ \ \ \texttt{Else} $L$ is not a
value set for any plane branch.\\ \hline
\end{tabular}
\end{center}

Let us illustrate the above algorithm with some examples.

\begin{example} Let us verify if the following are value set of $1$-forms for some plane branch:
$$\begin{array}{c}
L_1=\{6,9,12,15,16,17,18,21,22,24,25,27+\mathbb{N}\}, \\
L_2=\{6,9,12,15,16,17,18,21,22,23,24,25,27+\mathbb{N}\}, \\
L_3=\{6,9,12,15,16,18,19,21,22,23,24,25,27+\mathbb{N}\}, \\
L_4=\{6,9,12,15,16,18,19,21,22,24,25,27+\mathbb{N}\}.
\end{array}$$

The respective Ap\'ery set are:
$Ap(L_1)=Ap(L_2)=\{6,9,16,17,25,32\}$, $$Ap(L_3)=\{6,9,16,19,23,32\}\ \ \ \ \mbox{and}\ \ \ \  Ap(L_4)=\{6,9,16,19,29,32\}.$$

The sets $L_2, L_3$ and $L_4$ are covered by their Ap\'ery set, but $L_1$ is not covered by $Ap(L_1)$, because $23\in 17+\mathbb{N}\cdot 6$ but $23\not\in L_1$, so $L_1$ is not a value set of $1$-forms for any plane branch.

We get $\varepsilon_0=6$, $\varepsilon_1=3$, $\varepsilon_2=1$, $\eta_0=1$, $\eta_1=2$, $\eta_2=3$ for $L_i$ with $i=2,3,4$,
$$B_0(L_2)=B_0(L_3)=B_0(L_4)=\{6\},\ \ \ B_{1}(L_2)=B_{1}(L_3)=B_{1}(L_4)=\{9\},$$
$$B_2(L_2)=\{16,17\},\ \ \ B_2(L_3)=B_2(L_4)=\{16,19\}.$$
The condition given in Theorem \ref{bresinsky} is fulfill for $L_3$ and $L_4$, but not for $L_2$ because $17=\max(B_2(L_2))<\eta_1\cdot\max(B_1(L_2))=18$. Hence, $L_2$ is not a value set for any plane branch.

For $L_3$ and $L_4$ we get the value semigroup $\langle 6,9,19\rangle$ and by Example \ref{ex6919} we verify that $L_4$ is a value set for any plane branch but the same is false for $L_3$.
\end{example}

\vspace{0.35cm}

\begin{tabular}{lcl}
Abreu, M. O. R. & & Hernandes, M. E. \\
{\it osnar$@$outlook.com} & & {\it mehernandes$@$uem.br}\\
\end{tabular}
\vspace{0.5cm}

\hspace{1cm} Universidade Estadual de Maring\'a\vspace{0.25cm}

\hspace{1.5cm} Maring\'a - Paran\'a - Brazil

\end{document}